\documentclass{article}
\usepackage{arxiv}
\usepackage{comment}
\usepackage[utf8]{inputenc} % allow utf-8 input
\usepackage[T1]{fontenc}    % use 8-bit T1 fonts
\usepackage{hyperref}       % hyperlinks
\usepackage{url}            % simple URL typesetting
\usepackage{booktabs}       % professional-quality tables
\usepackage{amsmath}
\usepackage{amsfonts}       % blackboard math symbols
\usepackage{nicefrac}       % compact symbols for 1/2, etc.
\usepackage{microtype}      % microtypography
\usepackage{xcolor}
\usepackage{graphicx}
\usepackage{stmaryrd}
\usepackage{bm}
\usepackage[normalem]{ulem}
\usepackage[tickmarkheight=0.1cm]{todonotes}
\usepackage[ruled,vlined]{algorithm2e}

\usepackage{lineno}

\hypersetup{
colorlinks=true, % make the links colored
linkcolor=blue, % color TOC links in blue
urlcolor=red, % color URLs in red
citecolor=magenta,
linktoc=all % 'all' will create links for everything in the TOC
}

\newcommand{\backassign}{=:}

\newcommand{\nocomma}{}

\newcommand{\tmop}[1]{\ensuremath{\operatorname{#1}}}

\newcommand{\tmtextit}[1]{\text{{\itshape{#1}}}}

\newcommand{\poly}{\mathbb{P}}

\newcommand{\ud}{\mathrm{d}}

\newcommand{\q}{\tilde{u}}
\newcommand{\fq}{\tilde{f}}
\newcommand{\Fq}{\tilde{F}}

\newtheorem{remark}{Remark}

\title{
Equivalence of ADER and Lax-Wendroff in \\
DG / FR framework for linear problems
}
\author{
Arpit~Babbar \\
Centre for Applicable Mathematics\\
Tata Institute of Fundamental Research\\
Bangalore -- 560065\\
\texttt{arpit@tifrbng.res.in} \\
\And
Praveen~Chandrashekar \\
Centre for Applicable Mathematics\\
Tata Institute of Fundamental Research\\
Bangalore -- 560065\\
\texttt{praveen@math.tifrbng.res.in}
}
\begin{document}
\maketitle
%-------------------------------------------------------------------------------------
\begin{abstract}
ADER (Arbitrary high order by DERivatives) and Lax-Wendroff (LW) schemes are
two high order single stage methods for solving time dependent partial
differential equations. ADER is based on solving a locally implicit equation to obtain a space-time predictor solution while LW is based on an explicit Taylor's expansion in time.
We cast the corrector step of ADER Discontinuous Galerkin (DG) scheme into
an equivalent quadrature free Flux Reconstruction (FR) framework and then
show that the obtained ADER-FR scheme is equivalent to the LWFR scheme with
D2 dissipation numerical flux for linear problems. This also implies that the two schemes have the same Fourier stability limit for time step size. The equivalence is verified by
numerical experiments.
\end{abstract}
% \keywords{Conservation laws \and hyperbolic PDE\and Lax-Wendroff \and Discontinuous Galerkin \and flux reconstruction\and ADER}
\section{Introduction}
Spectral element methods like Flux Reconstruction (FR)~{\cite{Huynh2007}} and
Discontinuous Galerkin (DG)~{\cite{cockburn2000}} have been very successful
for solving advection dominated equations~{\cite{witherden2014,Ranocha2022}}.
Under the appropriate choice of quadrature and solution points, the DG
spectral element method can be cast in the FR framework~{\cite{Huynh2007}}. In
comparison to traditional semi-discrete FR / DG schemes, which mainly use
Runge-Kutta time integration, there are single stage methods which are
fully-discrete and allow the achievement of a naturally cache-blocking and
communication-avoiding scheme, reducing the amount of necessary MPI
communication to a minimum~{\cite{Dumbser2018}}. There are two single stage
discretizations in FR / DG framework that we discuss in this work. The first
approach is the family of Arbitrary high order schemes using DERivatives
(ADER) initially introduced by the idea of a generalized Riemann
solver~{\cite{Toro2001}} but later extended to Finite Volume / DG framework to
obtain high order accuracy by using a predictor-corrector
approach~{\cite{Dumbser2008}}. The second are the high
order Lax-Wendroff schemes introduced in~{\cite{Qiu2003,Qiu2005b}} which are
based on performing a Taylor's expansion in time and using the PDE to replace
temporal derivatives with spatial. The local evolution in ADER schemes is
performed by solving an element local implicit equation while the LW scheme
uses an explicit Taylor's expansion. In this work, we prove, for linear problems, the
equivalence of the ADER-DG scheme introduced in~{\cite{Dumbser2008}} with
Lax-Wendroff FR (LWFR) using D2 dissipation numerical flux introduced
in~{\cite{oldpaper}; the key observation used is that the space
time predictor polynomial can be explicitly determined for linear problems. We remark
that there are some works where both these ideas are considered as types of
ADER schemes, although in this work, we refer to ADER schemes as those which
use a local implicit solver like in~{\cite{Dumbser2008}} while LW schemes as
those that use a local Taylor's expansion~\cite{Qiu2005b,oldpaper}.
The rest of this paper is organized as follows. In Section~\ref{sec:ader.dg},
we review the ADER-DG scheme of~{\cite{Dumbser2008}} for 1-D scalar
conservation laws and cast it in an FR framework for simplicity of the proof.
In Section~\ref{sec:linear.equivalence}, we show the equivalence of
ADER-FR scheme and the LWFR scheme with D2 dissipation flux
of~{\cite{oldpaper}} for linear problems. In Section~\ref{sec:num}, we verify the equivalence
numerically and draw conclusions from the work in Section~\ref{sec:con}.
\section{ADER Discontinuous Galerkin and Flux
Reconstruction}\label{sec:ader.dg}
The arguments in this work apply to linear conservation laws of any dimension
but for simplicity we restrict ourselves to 1-D linear scalar conservation law
\begin{equation}
u_t + f (u)_x = 0, \qquad f (u) = au,\quad a = \text{const.} \label{eq:scalar.con.law}
\end{equation}
where $u$ is some conserved quantity, together with some initial and boundary
conditions. In this work, we consider the ADER-DG framework
of~{\cite{Dumbser2008}}. We will divide the physical domain $\Omega$ into
disjoint elements $\Omega_e$, with $\Omega_e = [x_{e - 1 / 2}, x_{e + 1 / 2}]$
and $\Delta x_e = x_{e + 1 / 2} - x_{e - 1 / 2}$. The temporal discretization
is performed by denoting the $n^{\tmop{th}}$ time interval as $[t_n, t_{n +
1}]$ and $\Delta t_n = t_{n + 1} - t_n$. Let us map all spatial and temporal
elements to reference elements $\Omega_e \to [0, 1]$, $[t_n, t_{n + 1}] \to
[0, 1]$ by
\[
x \mapsto \xi = \frac{x - x_{e - 1 / 2}}{\Delta x_e}, \qquad t \mapsto \tau
= \frac{t - t_n}{\Delta t_n}
\]
Thus, $x, t$ are physical variables in space and time and $\xi, \tau$ are the respective reference variables. Inside each element, we approximate the
solution as $\poly_N$ functions which are polynomials of degree $N \ge 0$. For
this, choose $N + 1$ distinct nodes $\{ \xi_i \}_{i = 0}^N$ in $[0, 1]$ which
will be taken to be Gauss-Legendre (GL) or Gauss-Lobatto-Legendre (GLL) nodes,
and will also be referred to as \tmtextit{solution points}. There are
associated quadrature weights $w_j$ such that the quadrature rule is exact for
polynomials of degree up to $2 N + 1$ for GL points and upto degree $2 N - 1$
for GLL points. Note that the nodes and weights we use are with respect to the
interval $[0, 1]$ whereas they are usually defined for the interval $[- 1, +
1]$. For constructing the space-time predictor, we use the same solution
points in time. The numerical solution inside an element $\Omega_e$ at $t =
t^n$ is given by
\[ x \in \Omega_e : \qquad u_h^n (\xi) = \sum_{j = 0}^N u_j^e \ell_j (\xi) \]
where each $\ell_j$ is a Lagrange polynomial of degree $N$ in $[0, 1]$ defined to satisfy $\ell_j (\xi_i) = \delta_{i \nocomma j}$ for $0 \leq i, j \leq N$.
\paragraph{Predictor step.}
The predictor inside a space-time element is given by
\begin{equation}
(x, t) \in \Omega_e \times [t_n, t_{n + 1}] : \qquad \q_h (\xi, \tau) =
\sum_{j, k = 0}^N \q_{j, k}^e \ell_j (\xi) \ell_k (\tau)
\label{eq:predictor.defn}
\end{equation}
Within each element $\Omega_e$, we take a local space-time test function
$\ell_{j, k}$
\[ (x, t) \in \Omega_e \times [t_n, t_{n + 1}] : \qquad \ell_{j, k} (\xi,
\tau) = \ell_j (\xi) \ell_k (\tau) \]
To compute the cell-local predictor, we multiply the conservation
law~\eqref{eq:scalar.con.law} by $\ell_{j, k}$ and do an integration by parts
in time
\begin{equation}
- \int_{t^n}^{t^{n + 1}} \int_{\Omega_e} \q_h \partial_t \ell_{j, k}  \ud x
\ud t +
 \int_{\Omega_e} \q_h (\xi, 1) \ell_{j, k}  \ud x -
\int_{\Omega_e} u_h^n (\xi) \ell_{j, k}  \ud x
+ \int_{t^n}^{t^{n +
1}} \int_{\Omega_e} (\partial_x  \fq_h) \ell_{j, k}  \ud x \ud t = 0
\label{eq:predictor.eqn}
\end{equation}
where $\fq_h = a \q_h$. The above system of equations~\eqref{eq:predictor.eqn}
is solved for all $\q_{j,k}^e$~\eqref{eq:predictor.defn}.
\paragraph{Corrector step.}
Integrate~\eqref{eq:scalar.con.law} over the space-time element $\Omega_e \times [t_n, t_{n+1}]$ with the test function $\ell_j = \ell_j (\xi)$ and perform an integration by parts in space to get
\begin{equation}
\begin{split}
\int_{\Omega_e} u_h^{n + 1} \ell_j  \ud x = \int_{\Omega_e} u_h^n \ell_j
\ud x & + \int_{t^n}^{t^{n + 1}} \int_{\Omega_e} \fq_h \partial_x \ell_j  \ud
x \ud t \\
& - \ell_j (1)  \int_{t^n}^{t^{n + 1}} f_{e + 1 / 2} (\q_h (t))
\ud t + \ell_j (0)  \int_{t^n}^{t^{n + 1}} f_{e - 1 / 2} \left( \q_h (t) \right)  \ud t \label{eq:ader.corr1}
\end{split}
\end{equation}
where, for the linear case, $f_{e + 1 / 2} \left( \q_h (t) \right)$ is the upwind flux
\begin{equation}
f_{e + 1 / 2} \left( \q_h (t) \right) = \frac{a}{2}  \left( \q_h
(x_{e + 1 / 2}^-, t) + \q_h (x_{e + 1 / 2}^+, t) \right) - \frac{| a |}{2}
\left( \q_h (x_{e + 1 / 2}^+, t) - \q_h (x_{e + 1 / 2}^-, t) \right)
\label{eq:numflux.defn}
\end{equation}
The complete numerical scheme is given by space-time quadrature on~\eqref{eq:ader.corr1}
at the solution points. By linearity of
the flux, quadrature on the flux term in~\eqref{eq:ader.corr1} is exact as we use GL / GLL quadrature points
and can thus perform another integration by parts in space to write
\begin{align*}
\int_{\Omega_e} u_h^{n + 1} \ell_j  \ud x = & \int_{\Omega_e} u_h^n \ell_j
\ud x - \int_{t^n}^{t^{n + 1}} \int_{\Omega_e} (\partial_x  \fq_h) \ell_j  \ud
x \ud t\\
& - \ell_j (1)  \int_{t^n}^{t^{n + 1}} (f_{e + 1 / 2} (\q_h (t)) -
\fq_h (1, t)) \ud t + \ell_j (0)  \int_{t^n}^{t^{n + 1}} (f_{e - 1 / 2}
(\q_h (t)) - \fq_h (0, t)) \ud t
\end{align*}
Performing quadrature in space at solution points gives
\begin{align*}
u_j^{n + 1} = & u_j^n -  \left[ \partial_x  \int_{t^n}^{t^{n + 1}}
\fq_h (t) \ud t \right] _j\\
& - \frac{\ell_j (1)}{\Delta x_e w_j}  \int_{t^n}^{t^{n + 1}} \left( f_{e + 1 / 2}
\left( \q_h (t) \right) - \fq_h (1, t) \right)  \ud t + \frac{\ell_j
(0)}{\Delta x_e w_j}  \int_{t^n}^{t^{n + 1}} \left( f_{e - 1 / 2} \left( \q_h (t) \right) - \fq_h (0, t) \right)  \ud t
\end{align*}
We choose correction functions $g_L, g_R \in \mathbb{P}_{N + 1}$ to be
$g_{\tmop{Radau}}$~{\cite{Huynh2007}} if the solution points are GL points and
$g_2$~{\cite{Huynh2007}} if solution points are GLL. Then, following the proof
of equivalence of Flux Reconstruction (FR) and DG in~{\cite{Huynh2007}}\footnote{For
Radau correction functions, the identities follow from taking $\phi = \ell_j$
in (7.8) of~{\cite{Huynh2007}}. For $g_2$, see Section 2.4
of~{\cite{Grazia2014}}.}
\[
g_R'(\xi_j) = \ell_j (1) / w_j, \qquad g_L' (\xi_j) = - \ell_j (0) / w_j
\] 
and thus the correction step can be written in the FR
form as
\begin{equation}
u_j^{n + 1} = u_j^n - \Delta t_n \partial_x  \Fq_h (\xi_j)
\label{eq:ader.evolution}
\end{equation}
where we define
\begin{equation}
\Fq_h (\xi) = \frac{1}{\Delta t_n}  \int_{t^n}^{t^{n + 1}} \left[ \fq_h
(\xi, t) + g_R (\xi)  \left( f_{e + 1 / 2} \left( \q_h (t) \right) -
\fq_h (1, t) \right) + g_L (\xi)  \left( f_{e - 1 / 2} \left( \q_h (t) \right) - \fq_h (0, t) \right) \right]  \ud t \label{eq:ader.fr.flux}
\end{equation}
which is the ADER time-averaged flux corrected by FR. The $g_L, g_R$ satisfy $g_L (0) = g_R (1) = 1, g_L (1) = g_R (0)
= 0$ so that 
\[
\Fq_h (0) =
\frac{1}{\Delta t_n}  \int_{t^n}^{t^{n + 1}} f_{e - 1 / 2} \left( \q_h
(t) \right)  \ud t, \qquad \Fq_h (1) = \frac{1}{\Delta t_n}
\int_{t^n}^{t^{n + 1}} f_{e + 1 / 2} \left( \q_h (t) \right)  \ud t
\]
making $\Fq_h$ a globally continuous flux approximation. The
equations~(\ref{eq:predictor.eqn},~\ref{eq:ader.evolution},~\ref{eq:ader.fr.flux}) describe the
ADER-FR scheme.
\section{Equivalence}\label{sec:linear.equivalence}
Since $f (u) = au$ in~\eqref{eq:scalar.con.law}, the numerical flux function
is linear and thus the corrected ADER time-averaged flux
of~\eqref{eq:ader.fr.flux} can be written as
\begin{equation}
\begin{split}
\Fq_h (\xi) = \frac{1}{\Delta t_n}  \int_{t^n}^{t^{n + 1}} a \q_h
(\xi, t)  \ud t & + g_R (\xi)  \left[ f_{e + 1 / 2} \left( \frac{1}{\Delta
t_n}  \int_{t^n}^{t^{n + 1}} \q_h (t)  \ud t \right) -
\frac{1}{\Delta t_n}  \int_{t^n}^{t^{n + 1}} a \q_h (1, t) \ud t \right]\\
& + g_L (\xi)  \left[ f_{e - 1 / 2} \left( \frac{1}{\Delta t_n}
\int_{t^n}^{t^{n + 1}} \q_h (t)  \ud t \right) - \frac{1}{\Delta
t_n}  \int_{t^n}^{t^{n + 1}} a \q_h (0, t) \ud t \right]
\end{split} \label{eq:ader.corr.flux}
\end{equation}
We will prove equivalence assuming that both schemes have the same
solution at time $t = t^n$. Now, by~{\cite{oldpaper}} (Equation (5)
of~{\cite{oldpaper}}, in particular), Lax-Wendroff Flux Reconstruction
(LWFR) in an element can be written as
\begin{equation}
u_j^{n + 1} = u_j^n - \Delta t_n \partial_x F_h (\xi_j)
\label{eq:lwfr.evolution}
\end{equation}
where $F_h$ is the continuous LW time averaged flux
\begin{equation}
F_h(\xi) = F_h^{\delta}(\xi) + g_R(\xi)  [F_{e + 1 / 2} - F_h^{\delta} (1)] + g_L(\xi)  [F_{e -
1 / 2} - F_h^{\delta} (0)] \label{eq:first.time.average.flux}
\end{equation}
and $F_h^{\delta}$ is the discontinuous time averaged flux computed by the approximate Lax-Wendroff procedure of~\cite{Zorio2017,Burger2017}, which gives the following for linear flux
\begin{equation}
F_h^{\delta} = \sum_{k = 0}^N \frac{\Delta t^k}{(k + 1) !}
\partial_t^k f (u^n) = a \sum_{k = 0}^N \frac{\Delta t^k}{(k + 1) !}
\partial_t^k u^n = a \sum_{k = 0}^N \frac{(- a \Delta t)^k}{(k + 1) !}
\partial_x^k u^n_h \backassign aU_h^n \label{eq:Uh.defn}
\end{equation}
where $U_h^n$ is the approximate time averaged solution, and all spatial
derivatives are computed as local polynomial derivatives. The numerical flux
with D2 dissipation introduced in~{\cite{oldpaper}} is given by
\begin{equation}
\begin{split}
F_{e + 1 / 2} & = \frac{1}{2}  (F_h^{\delta} (x_{e + 1 / 2}^-) +
F_h^{\delta} (x_{e + 1 / 2}^+)) - \frac{| a |}{2}  (U_h^n (x_{e + 1 /
2}^+) - U_h^n (x_{e + 1 / 2}^-))\\
& = \frac{a}{2}  (U_h^n (x_{e + 1 / 2}^-) + U_h^n (x_{e + 1 / 2}^+)) -
\frac{| a |}{2}  (U_h^n (x_{e + 1 / 2}^+) - U_h^n (x_{e + 1 / 2}^-))  \\
&= f_{e + 1 / 2} (U_h^n)
\end{split} \label{eq:numflux.d2}
\end{equation}
where $f_{e + 1 / 2} (U_h^n)$ is as defined
in~\eqref{eq:numflux.defn}. Thus, the time averaged flux~\eqref{eq:first.time.average.flux} in LWFR scheme~\eqref{eq:lwfr.evolution} can be written as
\begin{equation}
F_h(\xi) = aU_h^n(\xi) + g_R(\xi)  [f_{e + 1 / 2} (U_h^n) - U_h^n (1)] + g_L(\xi)  [f_{e
- 1 / 2} (U_h^n) - U_h^n (0)] \label{eq:lw.corr.flux}
\end{equation}
\begin{remark}
The D1 dissipation numerical flux, as termed in~{\cite{oldpaper}}, was used
in earlier works like~{\cite{Qiu2005b}} and is given by
\begin{equation}
F_{e + 1 / 2} = \frac{1}{2}  (F_h^\delta (x_{e + 1 / 2}^-) + F_h^\delta (x_{e + 1 / 2}^+)) -
\frac{| a |}{2}  (u_h^n (x_{e + 1 / 2}^+) - u_h^n (x_{e + 1 / 2}^-))
\label{eq:numflux.d1}
\end{equation}
The D2 flux~\eqref{eq:numflux.d2}, introduced in~{\cite{oldpaper}}, enhances the Fourier CFL stability limit. The equivalence between LW and ADER only holds with the D2 dissipation.
\end{remark}
Looking at~(\ref{eq:ader.evolution},~\ref{eq:lwfr.evolution}), to prove the
claimed equivalence, we need to show that~\eqref{eq:ader.corr.flux}
and~\eqref{eq:lw.corr.flux} are the same, which will be true if we show that
the time averaged solution $U_h^n$ defined in~\eqref{eq:Uh.defn} is given by
\begin{equation}
U_h^n(\xi) = \frac{1}{\Delta t_n}  \int_{t^n}^{t^{n + 1}} \q_h (\xi, t)  \ud
t \label{eq:final.claim}
\end{equation}
For simplicity of explanation, extend the cell local polynomial $x \mapsto \sum_{j = 0}^N u^e_j \ell_j (\xi (x))$ as a polynomial in whole of $\mathbb{R}$, now denoted $u_e^n$. Then, defined
in physical coordinates, $(x, t) \mapsto u^n_e (x - a (t - t_n))$ is a degree
$N$ space-time polynomial which satisfies the predictor
equation~\eqref{eq:predictor.eqn} for $f (u) = au$. Since the predictor
equation has a unique solution~{\cite{Dumbser2009}}, the solution
of~\eqref{eq:predictor.eqn} is indeed given in physical coordinates as
\begin{equation}
\q_h (x, t) = u^n_e (x - a (t - t^n)), \qquad x \in \Omega_e
\label{eq:predictor.linear.solution}
\end{equation}
Thus, we have $\partial_t  \q_h = - a \partial_x  \q_h$ and $\q_h |_{t = t_n,
x \in \Omega_e} = u_e^n = u_h^n$ which we will now exploit to
obtain~\eqref{eq:final.claim}. Since $\q_h$ is a degree $N$ polynomial, its
Taylor's expansion gives
\begin{align*}
\q_h (\xi, t) & = \sum_{k = 0}^N \frac{(t - t^n)^k}{k!} \partial_t^k  \q_h
(\xi, t^n)\\
& = \sum_{k = 0}^N \frac{(- a (t - t^n))^k}{k!} \partial_x^k  \q_h (\xi,
t^n) = \sum_{k = 0}^N \frac{(- a (t - t^n))^k}{k!} \partial_x^k u_h^n \qquad
\eqref{eq:predictor.linear.solution}\\
\implies\frac{1}{\Delta t_n}  \int_{t^n}^{t^{n + 1}} \q_h (\xi, t)  \ud
t & = \sum_{k = 0}^N \frac{(- a \Delta t)^k}{(k + 1) !} \partial_x^k
u_h^n = U_h^n \qquad \eqref{eq:Uh.defn}
\end{align*}
Thus, we have obtained~\eqref{eq:final.claim} proving equivalence of the two
schemes.
\begin{remark}
The above steps are not valid for a non-linear flux because the identity $\q_t = -
\fq_x$ need not hold at $t = t^n$.\label{rmk:non.linear}
\end{remark}
\section{Numerical validation}\label{sec:num}
The ADER-FR scheme described in Section~\ref{sec:ader.dg} is implemented,
tested and validated for general scalar conservation laws like
Burgers' equations with smooth solutions. For numerical validation of
equivalence, the Lax-Wendroff scheme with D1,
D2~(\ref{eq:numflux.d1},~\ref{eq:numflux.d2}) dissipation (called LW-D1,
LW-D2) and ADER scheme are tested for scalar linear advection
equation~\eqref{eq:scalar.con.law} with $a = 5$ and wave packet initial
condition $u (x, 0) = e^{- 10 x^2} \sin (10 \pi x)$ on domain $[- 1, 1]$ with
periodic and Dirichlet boundary conditions for degrees $N=1,2,3$. The non-periodic boundaries for LWFR are treated as in~\cite{oldpaper}. The Radau correction
functions~{\cite{Huynh2007}} and Gauss-Legendre solution points are used in the results shown, although we have also tested other correction functions
and solution points where same behavior was seen. Each scheme uses the same
time step size, and is within the stability limit. The LW-D2~\eqref{eq:numflux.d2} and ADER schemes are found to match
to $O (10^{- 14})$ in $L^{\infty}$ norm, verifying equivalence. In
Figure~\ref{fig:error}, we show the $L^2$ error $\| u_h - u_{\tmop{exact}}
\|_2$ versus time plot for LW-D1,
LW-D2~(\ref{eq:numflux.d1},~\ref{eq:numflux.d2}) and the ADER scheme for periodic (Figure~\ref{fig:error}a,~b,~c) and non-periodic (Figure~\ref{fig:error}d,~e,~f) boundaries. Since the ADER and LW-D2 schemes are equivalent, we see their $L^2$ error curves overlap, while for D1 dissipation, we see differences of upto $O(10^{- 2})$ for periodic boundaries and $O(10^{-3})$ for non-periodic boundaries. Thus, equivalence holds precisely with the D2 dissipation. The code used to generate these results is available online at
\href{https://github.com/Arpit-Babbar/ADER_FR}{https://github.com/Arpit-Babbar/ADER\_FR}.
\begin{figure}
\centering
\begin{tabular}{ccc}
{\includegraphics[width=0.28\textwidth]{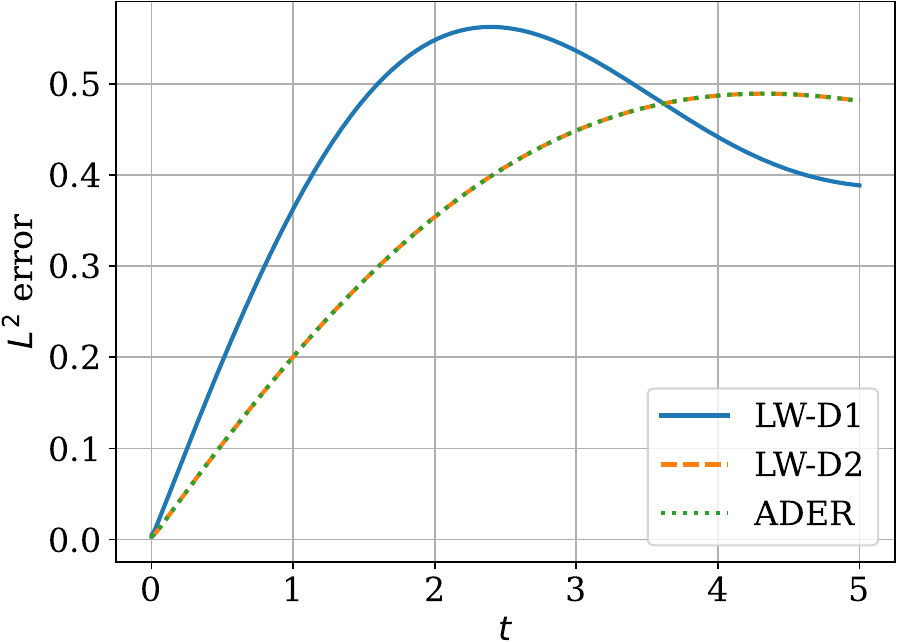}} &
{\includegraphics[width=0.28\textwidth]{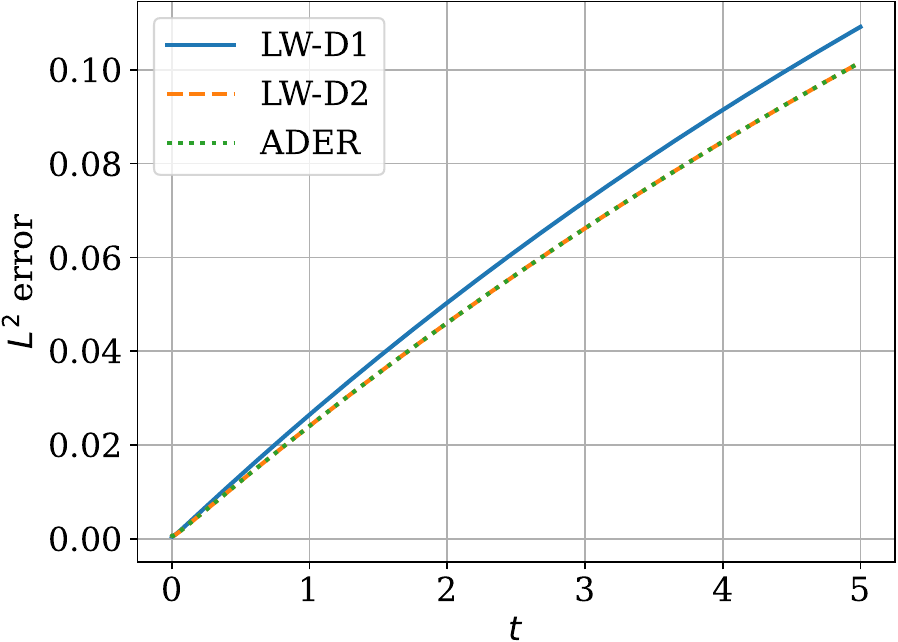}} &
{\includegraphics[width=0.28\textwidth]{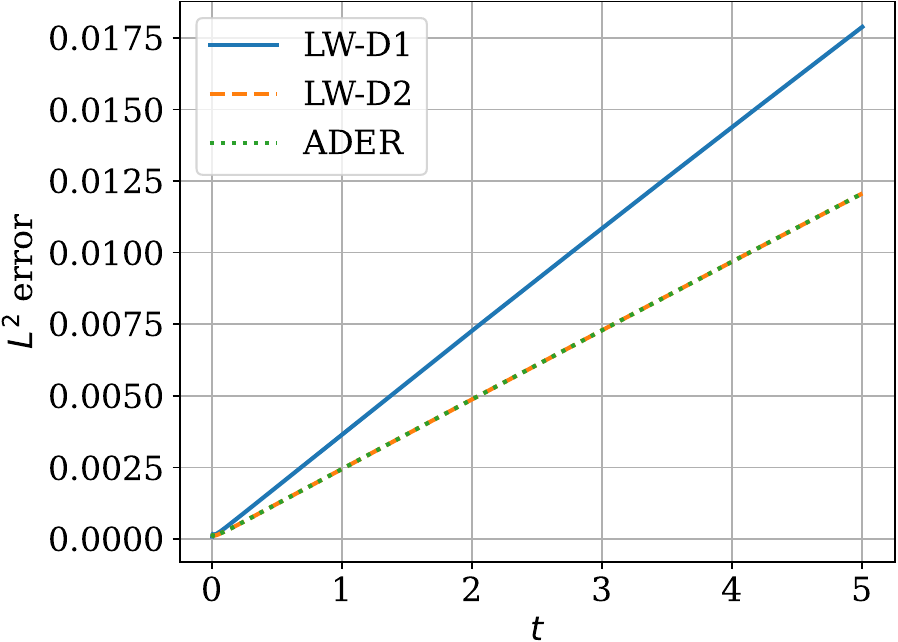}} \\
(a) & (b) & (c) \\
{\includegraphics[width=0.28\textwidth]{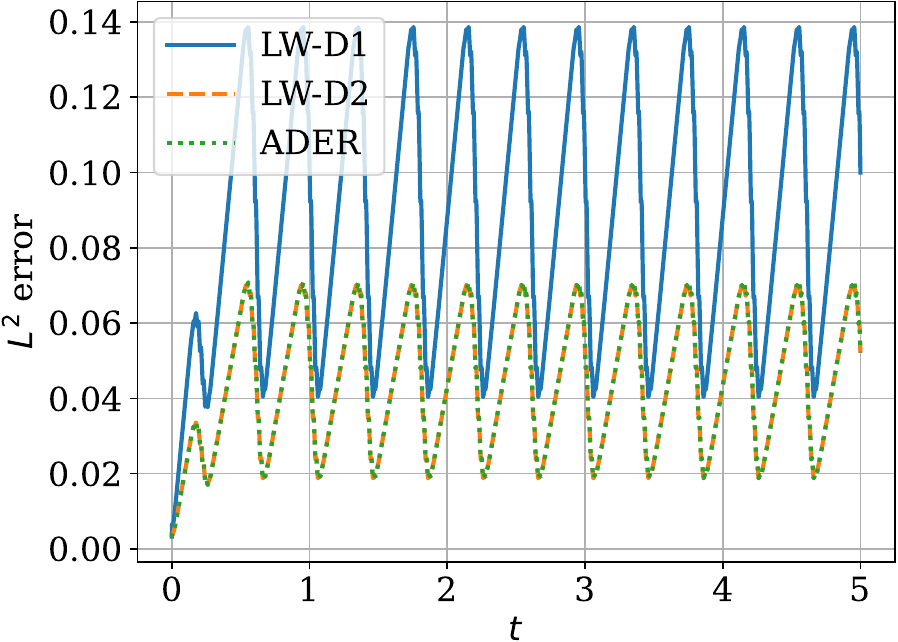}} &
{\includegraphics[width=0.28\textwidth]{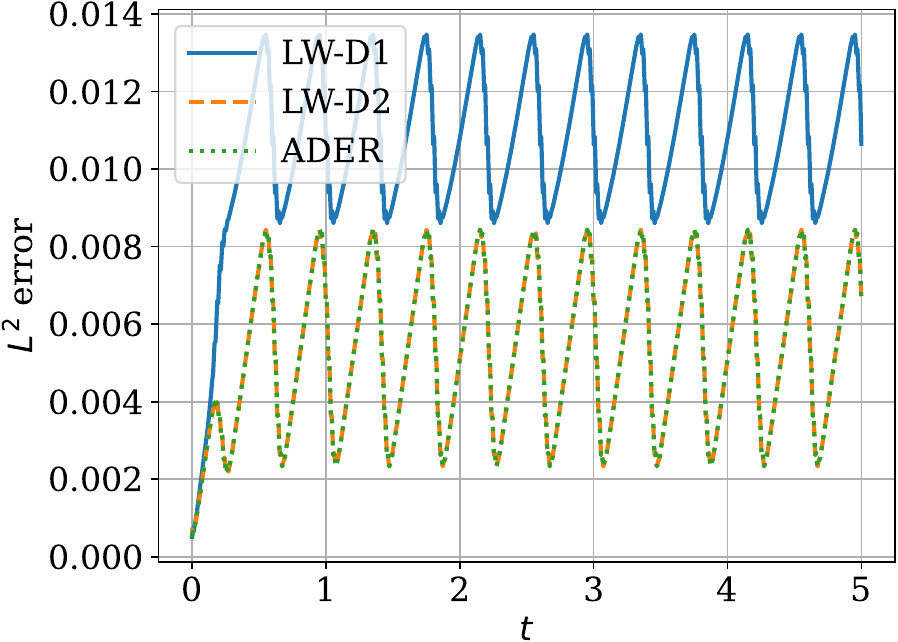}} &
{\includegraphics[width=0.28\textwidth]{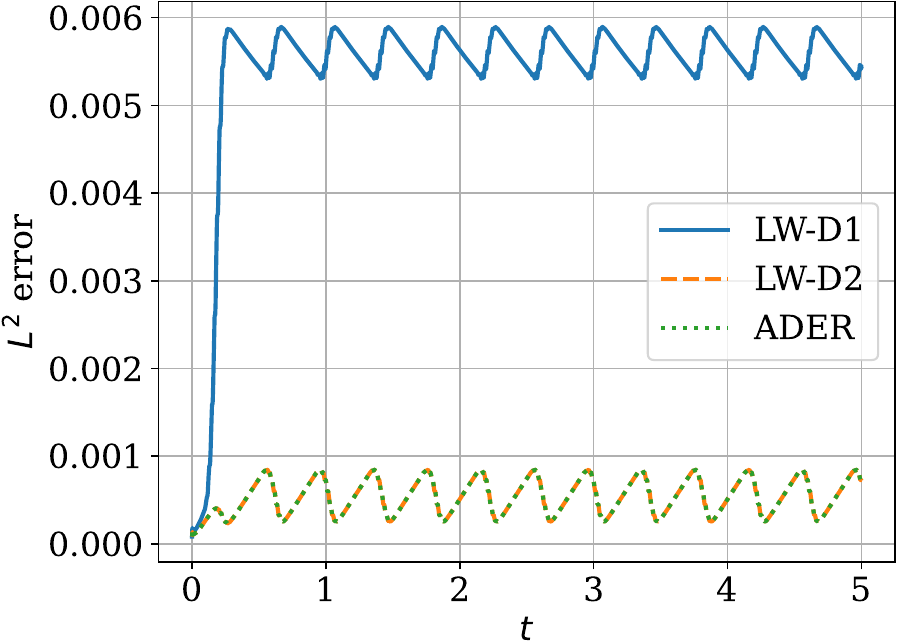}} \\
(d) & (e) & (f)
\end{tabular}
\caption{$L^2$ error $\| u_h - u_{\tmop{exact}} \|_2$ versus time for
wave packet test for different polynomial degrees with 240 degrees of freedom. Periodic : (a) $N = 1$, (b) $N = 2$,
(c) $N = 3$. Non-periodic : (a) $N=1$, (b) $N=2$, (c) $N=3$ \label{fig:error}}
\end{figure}
%\begin{figure}
%\centering
%\begin{tabular}{ccc}
%{\includegraphics[width=0.27\textwidth]{linear/energy1.pdf}} &
%{\includegraphics[width=0.27\textwidth]{linear/energy2.pdf}} &
%{\includegraphics[width=0.27\textwidth]{linear/energy3.pdf}} \\
%(a) & (b) & (c)
%\end{tabular}
%\caption{$L^2$ norm of numerical solution $u_h$ versus time for constant
%linear advection with periodic boundary conditions, initial condition $u (x,
%0) = e^{- 10 x^2} \sin (10 \pi x)$, $x \in [- 1, 1]$, for different
%polynomial degrees, each with 240 dofs. (a) $N = 1$, (b) $N = 2$, (c) $N =
%3$, (d) $N = 4$.\label{fig:energy}}
%\end{figure}
\section{Conclusions}\label{sec:con}
This work proves linear equivalence of high order ADER and Lax-Wendroff (LW)
schemes in Discontinuous Galerkin / Flux Reconstruction framework when the
numerical flux in LW is computed using the D2 dissipation introduced
in~{\cite{oldpaper}}. This is consistent with the Fourier stability analysis performed in~\cite{oldpaper} where it was observed that the CFL numbers of LWFR scheme with D2 dissipation are the same as those of ADER-DG schemes obtained in~\cite{Dumbser2008,Gassner2011a}. The equivalence was also numerically validated for a
wave packet test. The crucial observation needed for the proof is that the
solution of the predictor equation of ADER scheme has the same expression as
exact solution of the linear problem. Thus, this work relates two
single stage methods which are based on very different ideas and is thus a
contribution to our understanding of these numerical schemes. A natural but
important research question for further comparison of these two schemes is whether they
agree upto optimal order of accuracy, at least for smooth problems.
\bibliographystyle{siam}
\bibliography{references}

% \printbibliography
\end{document}